\documentclass[a4paper,10pt]{article}

\newtheorem{theorem}{Theorem}
\newtheorem{lemma}[theorem]{Lemma}

\usepackage{graphicx}

\usepackage{mathptmx}      

\usepackage{amssymb}
\usepackage{amsmath}
\usepackage[dvips]{epsfig}

\numberwithin{equation}{section}

\title{The Fourier Transforms of the Chebyshev and Legendre Polynomials}

\author{A S Fokas\\
        \small{School of Engineering and Applied Sciences, Harvard University,}\\
        \small{Cambridge, MA 02138, USA}\\
        \small{\emph{Permanent address:}$\vspace{-0.08cm}$}\\
        \small{Department of Applied Mathematics and Theoretical Physics,}\\
        \small{University of Cambridge, Cambridge CB3 0WA}\\
        \normalsize{$\vspace{-0.2cm}$}\\
        \large{S A Smitheman}\\
        \small{School of Biosciences, University of Nottingham, Sutton Bonington Campus,}\\
        \small{Loughborough, Leicestershire LE12 5RD, UK}
        }

\date{20 September 2012}

\begin{document}

\maketitle

\begin{abstract}
Analytic expressions for the Fourier transforms of the Chebyshev and Legendre polynomials are derived, and the latter is used to find a new \mbox{representation} for the half-order Bessel functions. The numerical \mbox{implementation} of the so-called unified method in the interior of a \mbox{convex} polygon provides an example of the applicability of these analytic expressions.
\end{abstract}

\section{Introduction}
The importance of Chebyshev and Legendre polynomials in numerical \mbox{analysis} is well known. The finite Fourier transform plays a crucial role in \mbox{mathematics} and applications. In spite of these facts it appears that the finite Fourier transform of Chebyshev and Legendre polynomials has not been constructed. Here, we present explicit formulae for these transforms.

The finite Fourier transform of the Legendre polynomials can be expressed in terms of the half-order Bessel functions, thus a direct application of our results is the construction of an explicit representation for the half-order Bessel functions.

Our results are motivated from the recent work of B Fornberg and collaborators \cite{D}, \cite{FF} on the numerical implementation of the unified method of \cite{F}-\cite{F3} to linear elliptic partial differential equations (PDEs) formulated in the interior of a convex polygon. In this case, the unified method yields a simple algebraic equation, the so-called global relation, which couples the finite Fourier transforms of the given boundary data with the finite Fourier transforms of the unknown boundary values. For the determination of these boundary values one has to choose (a) appropriate basis functions and (b) suitable collocation points in the Fourier space. Several such choices have appeared in the literature \cite{D}, \cite{FFSS}-\cite{SSF}; it appears that the best choice is (a) the unknown boundary values are expanded in terms of Legendre functions and (b) the collocation points are chosen on on rays introduced in \cite{FFSS}, \cite{SSF}. We note that Chebyshev and \mbox{Legendre} functions give rise to equivalent finite bases of $L^2$, and hence either choice will result in the same numerical method; however, in general the conditioning of the resulting linear systems will differ as conditioning is not invariant under matrix column operations.

As an example of the applicability of these analytic formulae we apply this approach to the simplest possible polygon, namely to a square, and we use the explicit construction of the finite Fourier transform of the Legendre functions instead of the usual representation in terms of the Bessel functions.

\section{Fourier transforms of the Chebyshev and \\Legendre polynomials}

\begin{theorem}
\label{Chebyshev}
Let $\hat T_m(\lambda)$ denote the finite Fourier transform of the Chebyshev polynomial $T_m(x)$, i.e.
\begin{equation}
\label{1C}
\hat T_m(\lambda)=\int_{-1}^1 e^{-i\lambda x}T_m(x)dx,\quad\lambda\in\mathbb C,\quad m=0,1,2,\ldots,
\end{equation}
where $T_m(x)$ denotes the Chebyshev polynomial
\begin{equation}
\label{2C}
T_m(x)=\cos(m\cos^{-1}(x)),\quad-1<x<1,\quad m=0,1,2,\ldots.
\end{equation}
Then,
\begin{equation}
\label{3C}
\hat T_m(0)=\left\{\begin{array}{ll}0,&m=1,\\\displaystyle\frac{(-1)^{m+1}-1}{m^2-1},&m=0,2,3,\ldots. \end{array}\right.
\end{equation}
Furthermore,
\begin{equation}
\label{4C}
\hat T_m(\lambda)=\sum_{n=1}^{m+1}\alpha_n^m\left[\frac{e^{i\lambda}}{(i\lambda)^n}+(-1)^{n+m}\frac{e^{-i\lambda}}{(i\lambda)^n}\right],\quad\lambda\in\mathbb C\smallsetminus\{0\},\quad m=0,1,2,\ldots,
\end{equation}
where the coefficients $\alpha_n^m$ are defined as follows:
\begin{equation}
\label{5C}
\alpha_1^m=(-1)^m,\quad\alpha_2^m=(-1)^{m+1}m^2,
\end{equation}
\begin{multline}
\label{6C}
\alpha_n^m=(-1)^{m+n-1}2^{n-2}m\sum_{k=1}^{m-n+2}\left(\!\begin{array}c n+k-3\\k-1\end{array}\!\right)\prod_{j=k}^{n+k-3}(m-j),\\n=3,4,\ldots,m+1.
\end{multline}
\end{theorem}
\newpage
\noindent\textbf{Proof} Making in \eqref{1C} the change of variables $x=\cos w$, we find
\begin{equation}
\label{6}
\hat T_m(\lambda)=\int_0^\pi e^{-i\lambda\cos w}\sin w\cos(mw)dw,\quad\lambda\in\mathbb C,\quad m=0,1,\ldots,
\end{equation}
from which the results for $\lambda=0$ follow. Using the identity
\begin{equation}
\label{7}
\sin we^{-i\lambda\cos w}=\frac1{i\lambda}\frac d{d w}\left(e^{-i\lambda\cos w}\right),
\end{equation}
and employing integration by parts in \eqref{6} we find
\begin{equation}
\label{8}
\hat T_m(\lambda)=\frac1{i\lambda}\left[e^{i\lambda}(-1)^m-e^{-i\lambda}+mK_m(-i\lambda)\right],\quad\lambda\in\mathbb C\smallsetminus\{0\},\quad m=0,1,\ldots,
\end{equation}
where
\begin{equation}
\label{9}
K_m(z)=\int_0^\pi e^{z\cos w}\sin(mw)dw,\quad z\in\mathbb C\smallsetminus\{0\},\quad m=0,1,\ldots.
\end{equation}
We now derive an analytic expression for $K_m(z)$. It is clear that $K_0(z)=0$, and on replacing $-i\lambda$ by $z$ in \eqref{7} we find
\begin{equation}
\label{10}
K_1(z)=\left.-\frac1ze^{z\cos w}\right|_{w=0}^\pi=\frac1z(e^z-e^{-z}),\quad z\in\mathbb C\smallsetminus\{0\}.
\end{equation}
The function $K_m(z)$ satisfies the relation
\begin{multline}
\label{11}
K_{m+1}(z)=-\frac{2m}zK_m(z)+K_{m-1}(z)+\frac2z(e^z+(-1)^{m-1})e^{-z}),\\
z\in\mathbb C\smallsetminus\{0\},\quad m=1,2,\ldots.
\end{multline}
For the derivation of this equation we will use the identity
\begin{equation}
\label{11a}
\sin((m+1)w)=\sin((m-1)w)+2\sin w\cos(mw).
\end{equation}
Thus,
\begin{equation}
\label{12}
\begin{split}
K_{m+1}(z)&=\int_0^\pi e^{z\cos w}\sin((m+1)w)dw\\
&=K_{m-1}(z)+2\int_0^\pi\sin we^{z\cos w}\cos(mw)dw.
\end{split}
\end{equation}
Using integration by parts we find that the above integral can be rewritten in the following form:
\begin{equation}
\label{13}
\begin{split}
&\int_0^\pi\sin we^{z\cos w}\cos(mw)dw\\
&\qquad=\left.-\frac1ze^{z\cos w}\cos(mw)\right|_0^\pi-\frac mz\int_0^\pi e^{z\cos w}\sin(mw)dw\\
&\qquad=\frac1z(e^z+(-1)^{m-1})e^{-z})-\frac mzK_m(z).
\end{split}
\end{equation}
We now show that for $m=0,1,\ldots$,
\begin{equation}
\label{14}
\begin{split}
&K_m(z)\\
&\quad=\frac mz(e^z+(-1)^me^{-z})\\
&\quad\phantom{=}+\frac1z\sum_{n=1}^{m-1}\left(\frac2z\right)^n((-1)^ne^z+(-1)^me^{-z})\sum_{k=1}^{m-n}\left(\!\begin{array}cn+k-1\\k-1\end{array}\!\right)\prod_{j=k}^{n+k-1}(m-j),
\end{split}
\end{equation}
from which \eqref{4C} follows using \eqref{8}. Clearly, this relation is valid for $m=0$ and $m=1$.

Furthermore, the function $K_m(z)$ in \eqref{14} satisfies the relation \eqref{11} if, and only if,
\begin{equation}
\begin{split}
\label{15}
&\frac1z\sum_{n=1}^m\left(\frac2z\right)^n((-1)^ne^z+(-1)^{m+1}e^{-z})\sum_{k=1}^{m-n+1}\left(\!\!\!\begin{array}cn+k-1\\k-1\end{array}\!\!\!\right)\prod_{j=k}^{n+k-1}(m+1-j)\\
&\,=-\frac{2m^2}{z^2}(e^z+(-1)^me^{-z})\\
&\,\phantom{=}-\frac{2m}{z^2}\sum_{n=1}^{m-1}\left(\frac2z\right)^{\!n}\!\!((-1)^ne^z+(-1)^me^{-z})\sum_{k=1}^{m-n}\left(\!\!\!\begin{array}cn+k-1\\k-1\end{array}\!\!\!\right)\prod_{j=k}^{n+k-1}(m-j)\\
&\,\phantom{=}+\frac1z\sum_{n=1}^{m-2}\!\left(\frac2z\right)^{\!n}\!\!\!((-1)^ne^z+(-1)^{m-1}e^{-z})\!\!\sum_{k=1}^{m-n-1}\!\!\left(\!\!\!\begin{array}cn+k-1\\k-1\end{array}\!\!\!\right)\!\!\prod_{j=k}^{n+k-1}\!\!(m-1-j).
\end{split}
\end{equation}
To prove \eqref{15}, we compare coefficients of each power of $z$ on each side of this equation. The coefficient of $z^{-2}$ on the left-hand side of \eqref{15} is given by
\begin{multline}
\label{16}
2(-e^z+(-1)^{m+1}e^{-z})\sum_{k=1}^mk(m+1-k)\\
=-(e^z+(-1)^me^{-z})\left[m(m+1)^2-\frac13m(m+1)(2m+1)\right],
\end{multline}
the coefficient of $z^{-2}$ on the right-hand side of \eqref{15} is given by
\begin{equation}
\label{17}
\begin{split}
&-(e^z+(-1)^me^{-z})\left[2m^2+2\sum_{k=1}^{m-2}k(m-1-k)\right]\\
&=-(e^z+(-1)^me^{-z})\!\!\left[2m^2+(m-1)^2(m-2)-\frac13(m-2)(m-1)(2m-3)\right]\!\!,
\end{split}
\end{equation}
and the expressions in \eqref{16}, \eqref{17} are equivalent. The coefficient of $z^{-(m+1)}$ on the left-hand side of \eqref{15} is given by
\begin{equation}
\label{18}
2^m((-1)^me^z+(-1)^{m+1}e^{-z})\prod_{j=1}^m(m+1-j)=2^m((-1)^me^z+(-1)^{m+1}e^{-z})m!,
\end{equation}
the coefficient of $z^{-(m+1)}$ on the right-hand side of \eqref{15} is given by
\begin{multline}
\label{19}
-2^mm((-1)^{m-1}e^z+(-1)^me^{-z})\prod_{j=1}^{m-1}(m-j)\\
=2^mm((-1)^me^z+(-1)^{m+1}e^{-z})(m-1)!,
\end{multline}
and the expressions in \eqref{18}, \eqref{19} are equivalent. The coefficient of $z^{-m}$ on the left-hand side of \eqref{15} is given by
\begin{equation}
\label{20}
\begin{split}
&2^{m-1}((-1)^{m-1}e^z+(-1)^{m+1}e^{-z})\sum_{k=1}^2\left(\!\begin{array}cm+k-2\\k-1\end{array}\!\right)\prod_{j=k}^{m+k-2}(m+1-j)\\
&\quad=2^{m-1}(-1)^{m-1}(e^z+e^{-z})\left[\prod_{j=1}^{m-1}(m+1-j)+m\prod_{j=2}^m(m+1-j)\right]\\
&\quad=2^mm!(-1)^{m-1}(e^z+e^{-z}),
\end{split}
\end{equation}
the coefficient of $z^{-m}$ on the right-hand side of \eqref{15} is given by
\begin{equation}
\label{21}
\begin{split}
&-2^{m-1}m((-1)^{m-2}e^z+(-1)^me^{-z})\sum_{k=1}^2\left(\!\begin{array}cm+k-3\\k-1\end{array}\!\right)\prod_{j=k}^{m+k-3}(m-j)\\
&\quad=2^{m-1}m(-1)^{m-1}(e^z+e^{-z})\left[\prod_{j=1}^{m-2}(m-j)+(m-1)\prod_{j=2}^{m-1}(m-j)\right]\\
&\quad=2^mm(-1)^{m-1}(e^z+e^{-z})(m-1)!,
\end{split}
\end{equation}
and the expressions in \eqref{20}, \eqref{21} are equivalent. Finally, we consider \\$z^{-(n+1)}$, $n=2,\,3,\,\ldots,\,m-2$. Omitting the common factor of \\$2^n((-1)^ne^z+(-1)^{m-1}e^{-z})$ for convenience, the difference of the coefficients of the two sides of equation \eqref{15} is given by
\begin{equation}
\label{22}
\begin{split}
&\sum_{k=1}^{m-n+1}\left[m\left(\!\begin{array}cn+k-2\\k-1\end{array}\!\right)\prod_{j=k}^{n+k-2}(m-j)-\left(\!\begin{array}cn+k-1\\k-1\end{array}\!\right)\prod_{j=k}^{n+k-1}(m+1-j)\right]\\
&\qquad\qquad\qquad\qquad\qquad\qquad\qquad+\sum_{k=1}^{m-n-1}\left(\!\begin{array}cn+k-1\\k-1\end{array}\!\right)\prod_{j=k}^{n+k-1}(m-1-j).
\end{split}
\end{equation}
Furthermore,
\begin{multline}
\label{23}
m\left(\!\begin{array}cn+k-2\\k-1\end{array}\!\right)\prod_{j=k}^{n+k-2}(m-j)-\left(\!\begin{array}cn+k-1\\k-1\end{array}\!\right)\prod_{j=k}^{n+k-1}(m+1-j)\\
=\left(\!\begin{array}cn+k-1\\k-1\end{array}\!\right)\left[\frac{mn}{n+k-1}-(m-k+1)\right]\prod_{j=k}^{n+k-2}(m-j),
\end{multline}
and
\begin{equation}
\label{24}
\begin{split}
&\sum_{k=1}^{m-n-1}\left(\!\begin{array}cn+k-1\\k-1\end{array}\!\right)\prod_{j=k}^{n+k-1}(m-1-j)\\
&\quad=\sum_{k=3}^{m-n+1}\left(\!\begin{array}cn+k-3\\k-3\end{array}\!\right)\prod_{j=k-2}^{n+k-3}(m-1-j)\\
&\quad=\sum_{k=1}^{m-n+1}\left(\!\begin{array}cn+k-1\\k-1\end{array}\!\right)\frac{(k-1)(k-2)(m-k+1)}{(n+k-1)(n+k-2)}\prod_{j=k}^{n+k-2}(m-j).
\end{split}
\end{equation}
Hence the above difference is given by
\begin{equation}
\label{25}
\begin{split}
&\sum_{k=1}^{m-n+1}\!\!\left[\rule{0cm}{0.7cm}\left(\!\begin{array}cn+k-1\\k-1\end{array}\!\right)\right.\\
&\phantom{\sum_{k=1}^{m-n+1}\!\!\left[\right.}\left.\left[\frac{mn}{n+k-1}-(m-k+1)\!\!\left(\!1-\frac{(k-1)(k-2)}{(n+k-1)(n+k-2)}\right)\!\right]\!\!\prod_{j=k}^{n+k-2}\!\!(m-j)\!\right]\\
&\,\,=n\sum_{k=1}^{m-n+1}\left[\rule{0cm}{0.7cm}\left(\!\begin{array}cn+k-1\\k-1\end{array}\!\right)\frac1{n+k-1}\right.\\
&\phantom{\,\,=n\sum_{k=1}^{m-n+1}\left[\right.}\left.\left[m-\frac{(n+2k-3)(m-k+1)}{n+k-2}\right]\prod_{j=k}^{n+k-2}(m-j)\right]\\
&\,\,=n\sum_{k=1}^{m-n}\left[\left(\!\begin{array}cn+k\\k\end{array}\!\right)\frac{k(n+2k-m-1)}{(n+k)(n+k-1)}\prod_{j=k+1}^{n+k-1}(m-j)\right]
=\frac1{(n-1)!}\sum_{k=1}^{m-n}a_k,
\end{split}
\end{equation}
where
\begin{equation}
\label{26}
a_k=\frac{(m-k-1)!(n+k-2)!(n+2k-m-1)}{(m-n-k)!(k-1)!},\quad k=1,\,2,\,\ldots,\,m-n.
\end{equation}
Hence
$a_{m-n+1-k}=-a_k$. If $m-n$ is odd, then
\begin{equation}
\label{27}
n+2\left(\frac{m-n+1}2\right)-m-1=0,
\end{equation}
and hence $a_{\frac{m-n+1}2}=0$. Thus the coefficients of $z^{-(n+1)}$ of the two sides of \eqref{15} are equal for $n=2,\,3,\,\ldots,\,m-2$. This establishes the validity of \eqref{15}, which completes the proof of \eqref{4C}.$\hfill\square$
\newpage
\begin{theorem}
\label{Legendre}
Let $\hat P_m(\lambda)$ denote the finite Fourier transform of the Legendre polynomial $P_m(x)$, i.e.
\begin{equation}
\label{1L}
\hat P_m(\lambda)=\int_{-1}^1 e^{-i\lambda x}P_m(x)dx,\quad\lambda\in\mathbb C,\quad m=0,1,2,\ldots,
\end{equation}
where $P_m(x)$ denotes the Legendre polynomial
\begin{equation}
\label{2L}
P_m(x)=\frac1{2^mm!}\frac{d^m}{dx^m}((x^2-1)^m),\quad-1<x<1,\quad m=0,1,2,\ldots.
\end{equation}
Then,
\begin{equation}
\label{3L}
\hat P_m(0)=\left\{\begin{array}{ll}2,&m=0,\\0,&m=1,2,\ldots. \end{array}\right.
\end{equation}
Furthermore,
\begin{equation}
\label{4L}
\hat P_m(\lambda)=\sum_{n=1}^{m+1}\beta_n^m\left[\frac{e^{i\lambda}}{(i\lambda)^n}+(-1)^{m+n}\frac{e^{-i\lambda}}{(i\lambda)^n}\right],\quad\lambda\in\mathbb C\smallsetminus\{0\},\quad m=0,1,2,\ldots,
\end{equation}
where the coefficients $\beta_n^m$ are defined as follows:
\begin{subequations}
\label{betas}
\begin{equation}
\label{5L}
\beta_1^m=1\textrm{ if $m$ even},
\end{equation}
\begin{equation}
\label{6L}
\begin{split}
&\beta_n^m=\left((m+n)\left(\!\begin{array}c\frac{m+n-3}2\\n-1\end{array}\!\right)+\left(\!\begin{array}c\frac{m+n-3}2\\n-2\end{array}\!\right)\right)R_{\frac{m-n+3}2,\frac{m+n-3}2}(m),\\
&\qquad\qquad\qquad\qquad\qquad\qquad\qquad\qquad n=\left\{\begin{array}{ll}3,5,...,m+1&\!\!\!\!\textrm{if $m$ even},\\2,4,...,m+1&\!\!\!\!\textrm{if $m$ odd},\end{array}\right.
\end{split}
\end{equation}
\begin{equation}
\label{7L}
\beta_n^m=-\left(\!\begin{array}c\frac{m+n}2-1\\n-1\end{array}\!\right)R_{\frac{m-n}2+1,\frac{m+n}2-1}(m),\,\, n=\left\{\begin{array}{ll}2,4,...,m&\!\!\!\!\textrm{if $m$ even},\\1,3,...,m&\!\!\!\!\textrm{if $m$ odd},\end{array}\right.
\end{equation}
\end{subequations}
where
\begin{equation}
\label{5b}
R_{s,t}(r)=\left\{\begin{array}{ll}\displaystyle\prod_{k=s}^t(2(r-k)+1),&s\leq t,\\1,&s>t,\end{array}\right.
\,\textrm{and}\,\left(\!\begin{array}c s\\t\end{array}\!\right)=\left\{\begin{array}{ll}\displaystyle\frac{s!}{(s-t)!t!},&s\geq t,\\0,&s<t.\end{array}\right.
\end{equation}
\end{theorem}
\textbf{Proof} We first consider \eqref{3L} for $m=0$:
\begin{equation}
\label{28}
\hat P_0(0)=\int_{-1}^1 P_0(x)dx=2.
\end{equation}
Using the identity
\begin{equation}
\label{29}
(2m+1)P_m(x)=\frac d{dx}(P_{m+1}(x)-P_{m-1}(x)),\quad-1<x<1,\quad m=1,2,\ldots,
\end{equation}
as well as the equations $P_m(1)=1$ and $P_m(-1)=(-1)^m$ for $m=0,1,\ldots$ we find the second line of \eqref{3L}:
\begin{equation}
\label{30}
\hat P_m(0)
=\int_{-1}^1P_m(x)dx=\frac1{2m+1}[P_{m+1}(x)-P_{m-1}(x)]_{-1}^1=0.
\end{equation}
For $m=0,\,1$, $P_m(\lambda)=T_m(\lambda)$ and thus $\hat P_m(\lambda)=\hat T_m(\lambda)$; hence \eqref{4C} establishes the validity of \eqref{4L}.

Furthermore, using \eqref{29} and integrating by parts we find
\begin{equation}
\label{31}
\hat P_{m+1}(\lambda)=-\frac i\lambda(2m+1)\hat P_m(\lambda)+\hat P_{m-1}(\lambda),\quad
\lambda\in\mathbb C\smallsetminus\{0\},\quad m=1,2,\ldots.
\end{equation}
For convenience, we consider the following equivalent forms of \eqref{4L} for $m$ even and odd respectively:
\begin{equation}
\label{32}
\begin{split}
\hat P_m(\lambda)&=\frac{(-1)^{\frac m2}}\lambda\sum_{j=1}^{\frac m2-1}\frac{(-1)^j}{\lambda^{m-2j-1}}\left((2(m-j)+1)\left(\!\begin{array}c m-j-1\\m-2j\end{array}\!\right)\psi(\lambda)\right.\\
&\phantom{=\frac{(-1)^{\frac m2}}\lambda\sum_{j=1}^{\frac m2-1}\frac{(-1)^j}{\lambda^{m-2j-1}}\quad}\left.+\left(\!\begin{array}c m-j-1\\m-2j-1\end{array}\!\right)\chi(\lambda)\right)R_{j+1,m-j-1}(m)\\
&\phantom{=}+\psi(\lambda)+\frac{(-1)^{\frac m2}}{\lambda^m}\chi(\lambda)R_{1,m-1}(m),
\end{split}
\end{equation}
and
\begin{equation}
\label{33}
\begin{split}
\hat P_m(\lambda)&=-\frac i\lambda\left(1+\frac{(-1)^{\frac{m-1}2}}{\lambda^{m-1}}R_{1,m-1}(m)\right)\chi(\lambda)-\frac{i(m-1)(m+2)}{2\lambda}\psi(\lambda)\\
&\phantom{=}-\frac{(-1)^{\frac{m-1}2}i}\lambda\sum_{j=1}^{\frac{m-3}2}\frac{(-1)^j}{\lambda^{m-2j-1}}\left((2(m-j)+1)\left(\!\begin{array}c m-j-1\\m-2j\end{array}\!\right)\psi(\lambda)\right.\\
&\phantom{=-\frac{(-1)^{\frac{m-1}2}i}\lambda\sum_{j=1}^{\frac{m-3}2}\frac{(-1)^j}{\lambda^{m-2j-1}}\left.\right.\,\,}\left.+\left(\!\begin{array}c m-j-1\\m-2j-1\end{array}\!\right)\chi(\lambda)\!\!\right)R_{j+1,m-j-1}(m),
\end{split}
\end{equation}
where
\begin{equation}
\label{34}
\psi(\lambda)=\frac i\lambda(e^{-i\lambda}-e^{i\lambda})\quad\textrm{and}\quad\chi(\lambda)=\psi(\lambda)-(e^{-i\lambda}+e^{i\lambda}).
\end{equation}
For odd $m$, \eqref{31} is equivalent to the following expression:
\begin{equation}
\label{35}
\begin{array}l
\displaystyle\!\frac{(-1)^{\frac{m+1}2}}\lambda\sum_{j=1}^{\frac{m-1}2}\frac{(-1)^j}{\lambda^{m-2j}}R_{j+1,m-j}(m+1)\vspace{-0.2cm}\\
\displaystyle\!\phantom{\frac{(-1)^{\frac{m+1}2}}\lambda\sum_{j=1}^{\frac{m-1}2}}\quad\left(\!\!(2(m-j)+3)\left(\!\!\!\!\begin{array}c m-j\\m-2j+1\end{array}\!\!\!\!\right)\psi(\lambda)+\left(\!\!\!\!\begin{array}c m-j\\m-2j\end{array}\!\!\!\!\right)\chi(\lambda)\!\!\right)\vspace{0cm}\\
\displaystyle\!+\psi(\lambda)+\frac{(-1)^{\frac{m+1}2}}{\lambda^{m+1}}R_{1,m}(m+1)\chi(\lambda)\vspace{0.2cm}\\
\displaystyle\!=\frac{(-1)^{\frac{m-1}2}}\lambda\!\sum_{j=1}^{\frac{m-3}2}\frac{(-1)^j}{\lambda^{m-2j-2}}R_{j+1,m-j-2}(m-1)\vspace{-0.2cm}\\
\displaystyle\!\phantom{=\frac{(-1)^{\frac{m-1}2}}\lambda\!\sum_{j=1}^{\frac{m-3}2}}
\left(\!\!(2(m-j)-1)\!\left(\!\!\!\!\begin{array}c m-j-2\\m-2j-1\end{array}\!\!\!\!\right)\!\psi(\lambda)+\!\left(\!\!\!\!\begin{array}c m-j-2\\m-2j-2\end{array}\!\!\!\!\right)\!\chi(\lambda)\!\!\right)\vspace{0.2cm}\\
\displaystyle\!\phantom{=}+\psi(\lambda)+\frac{(-1)^{\frac{m-1}2}}{\lambda^{m-1}}R_{1,m-2}(m-1)\chi(\lambda)\vspace{0.2cm}\\
\displaystyle\!\phantom{=}-\frac1{\lambda^2}(2m+1)\left(\!\!1+\frac{(-1)^{\frac{m-1}2}}{\lambda^{m-1}}R_{1,m-1}(m)\!\!\right)\chi(\lambda)\vspace{0.2cm}\\
\displaystyle\!\phantom{=}-\frac1{2\lambda^2}(2m+1)(m-1)(m+2)\psi(\lambda)\vspace{0.2cm}\\
\displaystyle\!\phantom{=}-\frac{(-1)^{\frac{m-1}2}}{\lambda^2}(2m+1)\vspace{0.2cm}\\
\displaystyle\!\phantom{=-}\quad\cdot\sum_{j=1}^{\frac{m-3}2}\frac{(-1)^j}{\lambda^{m-2j-1}}R_{j+1,m-j-1}(m)\vspace{-0.2cm}\\
\displaystyle\!\phantom{=-\quad\cdot\sum_{j=1}^{\frac{m-3}2}\quad}\left(\!\!\left(\!\!\!\!\begin{array}c m-j-1\\m-2j-1\end{array}\!\!\!\!\right)\chi(\lambda)+(2(m-j)+1)\left(\!\!\!\!\begin{array}c m-j-1\\m-2j\end{array}\!\!\!\!\right)\psi(\lambda)\!\!\right).
\end{array}
\end{equation}
Furthermore,
\begin{equation}
\label{36}
\begin{array}l
\displaystyle\sum_{j=1}^{\frac{m-3}2}\frac{(-1)^j}{\lambda^{m-2j-2}}\left((2(m-j)-1)\left(\!\begin{array}c m-j-2\\m-2j-1\end{array}\!\right)\psi(\lambda)\right.\vspace{-0.2cm}\\
\displaystyle\phantom{\sum_{j=1}^{\frac{m-3}2}\frac{(-1)^j}{\lambda^{m-2j-2}}\left.\right.}\quad\left.+\left(\!\begin{array}c m-j-2\\m-2j-2\end{array}\!\right)\chi(\lambda)\right)R_{j+1,m-j-2}(m-1)\vspace{-0.2cm}\\
\displaystyle\quad=\sum_{j=2}^{\frac{m-1}2}\frac{(-1)^{j-1}}{\lambda^{m-2j}}\left((2(m-j)+1)\left(\!\begin{array}c m-j-1\\m-2j+1\end{array}\!\right)\psi(\lambda)\right.\vspace{-0.2cm}\\
\displaystyle\phantom{\quad=\sum_{j=2}^{\frac{m-1}2}\frac{(-1)^{j-1}}{\lambda^{m-2j}}\left.\right.}\quad\left.+\left(\!\begin{array}c m-j-1\\m-2j\end{array}\!\right)\chi(\lambda)\right)R_{j,m-j-1}(m-1).
\end{array}
\end{equation}
\newpage
\noindent It is straightforward to show that
\begin{multline}
\label{37}
R_{j+1,m-j}(m+1)=(2(m-j)+1)R_{j+1,m-j-1}(m),\\R_{j,m-j-1}(m-1)=(2j+1)R_{j+1,m-j-1}(m).
\end{multline}
Hence for odd $m$, \eqref{31} is equivalent to the following equation:
\begin{equation}
\label{38}
\begin{array}l
\displaystyle\!\frac{(-1)^{\frac{m+1}2}}\lambda\sum_{j=1}^{\frac{m-1}2}\frac{(-1)^j}{\lambda^{m-2j}}(2(m-j)+1)R_{j+1,m-j-1}(m)\vspace{-0.2cm}\\
\displaystyle\!\phantom{\frac{(-1)^{\frac{m+1}2}}\lambda\sum_{j=1}^{\frac{m-1}2}}\,\,\left(\!\!(2(m-j)+3)\left(\!\!\!\!\begin{array}c m-j\\m-2j+1\end{array}\!\!\!\!\right)\!\psi(\lambda)+\left(\!\!\!\!\begin{array}c m-j\\m-2j\end{array}\!\!\!\!\right)\!\chi(\lambda)\!\!\right)\vspace{0cm}\\
\displaystyle\!+\frac{(-1)^{\frac{m+1}2}}{\lambda^{m+1}}(2m+1)R_{1,m-1}(m)\chi(\lambda)\vspace{0.2cm}\\
\displaystyle\!=\frac{(-1)^{\frac{m+1}2}}\lambda\sum_{j=2}^{\frac{m-1}2}\frac{(-1)^j}{\lambda^{m-2j}}(2j+1)R_{j+1,m-j-1}(m)\vspace{-0.2cm}\\
\displaystyle\!\phantom{=\frac{(-1)^{\frac{m+1}2}}\lambda\sum_{j=2}^{\frac{m-1}2}}\,\,\left(\!\!(2(m-j)+1)\left(\!\!\!\!\begin{array}c m-j-1\\m-2j+1\end{array}\!\!\!\!\right)\!\psi(\lambda)+\left(\!\!\!\!\begin{array}c m-j-1\\m-2j\end{array}\!\!\!\!\right)\!\chi(\lambda)\!\!\right)\vspace{0.2cm}\\
\displaystyle\!\phantom{=}+\frac{3(-1)^{\frac{m-1}2}}{\lambda^{m-1}}R_{2,m-2}(m)\chi(\lambda)\vspace{0.2cm}\\
\displaystyle\!\phantom{=}-\frac1{\lambda^2}(2m+1)\!\left(\!\!1+\frac{(-1)^{\frac{m-1}2}}{\lambda^{m-1}}R_{1,m-1}(m)\!\!\right)\!\chi(\lambda)\vspace{0.2cm}\\
\displaystyle\!\phantom{=}-\frac1{2\lambda^2}(2m+1)(m-1)(m+2)\psi(\lambda)\vspace{0.2cm}\\
\displaystyle\!\phantom{=}+\frac{(-1)^{\frac{m+1}2}}\lambda(2m+1)\vspace{0.2cm}\\
\displaystyle\!\phantom{=+}\quad\cdot\sum_{j=1}^{\frac{m-3}2}\frac{(-1)^j}{\lambda^{m-2j}}R_{j+1,m-j-1}(m)\vspace{-0.2cm}\\
\displaystyle\!\phantom{=+\quad\cdot\sum_{j=1}^{\frac{m-3}2}}\,\,\left(\!\!\left(\!\!\!\!\begin{array}c m-j-1\\m-2j-1\end{array}\!\!\!\!\right)\!\chi(\lambda)+(2(m-j)+1)\left(\!\!\!\!\begin{array}c m-j-1\\m-2j\end{array}\!\!\!\!\right)\!\psi(\lambda)\!\!\right).
\end{array}
\end{equation}
Hence,
\begin{equation}
\label{39}
\begin{array}l
\displaystyle\!\frac{(-1)^{\frac{m+1}2}}\lambda\vspace{0.2cm}\\
\displaystyle\!\cdot\sum_{j=2}^{\frac{m-3}2}\frac{(-1)^j}{\lambda^{m-2j}}R_{j+1,m-j-1}(m)\\
\displaystyle\!\phantom{\cdot\sum_{j=2}^{\frac{m-3}2}}\!\!\bigg(\!\!(2(m-j)+1)\!\left(\!\!(2(m-j)+3)\!\!\left(\!\!\!\!\begin{array}c m-j\\m-2j+1\end{array}\!\!\!\!\right)\!-(2j+1)\!\!\left(\!\!\!\!\begin{array}c m-j-1\\m-2j+1\end{array}\!\!\!\!\right)\right.\\
\displaystyle\!\qquad\qquad\qquad\qquad\qquad\qquad\qquad\qquad\quad\left.-(2m+1)\!\left(\!\!\!\!\begin{array}c m-j-1\\m-2j\end{array}\!\!\!\!\right)\!\!\right)\!\psi(\lambda)\\
\displaystyle\!\phantom{\,\cdot\sum_{j=2}^{\frac{m-3}2}\bigg(\!\!}+\!\left(\!\!(2(m-j)+1)\!\left(\!\!\!\!\begin{array}c m-j\\m-2j\end{array}\!\!\!\!\right)\!-(2j+1)\!\left(\!\!\!\!\begin{array}c m-j-1\\m-2j\end{array}\!\!\!\!\right)\right.\\
\displaystyle\!\qquad\qquad\qquad\qquad\qquad\qquad\qquad\qquad\quad\!\!\left.-(2m+1)\!\left(\!\!\!\!\begin{array}c m-j-1\\m-2j-1\end{array}\!\!\!\!\right)\!\!\right)\!\chi(\lambda)\!\!\bigg)\\
\displaystyle\!+\frac{(-1)^{\frac{m-1}2}}{\lambda^{m-1}}R_{2,m-2}(m)[(2m-1)((2m+1)\psi(\lambda)+(m-1)\chi(\lambda))\\
\displaystyle\!\phantom{+\frac{(-1)^{\frac{m-1}2}}{\lambda^{m-1}}R_{2,m-2}(m)[}\qquad\quad-(2m+1)((m-2)\chi(\lambda)+(2m-1)\psi(\lambda))]\\
\displaystyle\!-\frac1{\lambda^2}\!\left(\!\!(m+2)\!\left(\!\frac18(m+4)(m+1)(m-1)\psi(\lambda)+\frac12(m+1)\chi(\lambda)\!\right)\right.\\
\displaystyle\!\qquad\qquad\qquad\qquad\!\left.-m\!\left(\!\frac18(m+2)(m-1)(m-3)\psi(\lambda)+\frac12(m-1)\chi(\lambda)\!\right)\!\!\right)\!\\
\displaystyle\!+\!\left(\!\!\frac{(-1)^{\frac{m+1}2}}{\lambda^{m+1}}(2m+1)R_{1,m-1}(m)+\frac{3(-1)^{\frac{m+1}2}}{\lambda^{m-1}}R_{2,m-2}(m)\right.\\
\displaystyle\!\qquad\qquad\qquad\qquad\qquad\quad\left.+\frac1{\lambda^2}(2m+1)\!\left(\!\!1+\frac{(-1)^{\frac{m-1}2}}{\lambda^{m-1}}R_{1,m-1}(m)\!\!\right)\!\!\right)\!\chi(\lambda)\\
\displaystyle\!+\frac1{2\lambda^2}(2m+1)(m-1)(m+2)\psi(\lambda)\vspace{0.2cm}\\
\displaystyle\!=0.
\end{array}
\end{equation}
Equation \eqref{39} follows from the identities
\begin{multline}
\label{40}
(2(m-j)+3)\left(\!\!\!\begin{array}c m-j\\m-2j+1\end{array}\!\!\!\right)-(2j+1)\left(\!\!\!\begin{array}c m-j-1\\m-2j+1\end{array}\!\!\!\right)\\
-(2m+1)\left(\!\!\!\begin{array}c m-j-1\\m-2j\end{array}\!\!\!\right)=0,
\end{multline}
\begin{multline}
\label{41}
(2(m-j)+1)\left(\!\!\!\begin{array}c m-j\\m-2j\end{array}\!\!\!\right)-(2j+1)\left(\!\!\!\begin{array}c m-j-1\\m-2j\end{array}\!\!\!\right)\\
-(2m+1)\left(\!\!\!\begin{array}c m-j-1\\m-2j-1\end{array}\!\!\!\right)=0,
\end{multline}
\newpage
\noindent and the coefficients of $\displaystyle\frac{(-1)^{\frac{m+1}2}}{\lambda^{m+1}}R_{1,m-1}(m)\chi(\lambda)$, $\displaystyle\frac{(-1)^{\frac{m-1}2}}{\lambda^{m-1}}\psi(\lambda)R_{2,m-2}(m)$, \\$\displaystyle\frac{(-1)^{\frac{m-1}2}}{\lambda^{m-1}}\chi(\lambda)R_{2,m-2}(m)$, $-\displaystyle\frac1{8\lambda^2}\psi(\lambda)$ and $-\displaystyle\frac1{2\lambda^2}\chi(\lambda)$ on the left hand side of \eqref{39} vanishing:
\begin{gather*}
(2m+1)-(2m+1)=0,\,\,\,\,(2m-1)(2m+1)-(2m+1)(2m-1)=0,\\
(2m-1)(m-1)-(2m+1)(m-2)-3=0,\\
(m-1)(m+2)((m+1)(m+4)-m(m-3)-4(2m+1))=0,\\
(m+1)(m+2)-m(m-1)-2(2m+1)=0.
\end{gather*}
The validity of \eqref{31} for even $m$ can be demonstrated similarly, completing the proof of equation \eqref{4L}.$\hfill\square$

\section{An explicit representation for the half-order\\Bessel functions}

\begin{lemma} Let $J_{m+\frac12}(\lambda)$ denote the half-order Bessel function, i.e.
\begin{equation}
\begin{split}
J_{m+\frac12}(\lambda)&=\frac1\pi\int_0^\pi\cos\left(\left(m+\frac12\right)\tau-\lambda\sin\tau\right)d\tau\\
&+\frac{(-1)^{m+1}}\pi\int_0^\infty\exp\left(-\lambda\sinh\tau-\left(m+\frac12\right)\tau\right)d\tau,\\
&\qquad\qquad\qquad\qquad\qquad\qquad\qquad\lambda\in\mathbb C,\quad m=0,1,2,\ldots.
\end{split}
\end{equation}
Then
\begin{equation}
J_{m+\frac12}(0)=0,\quad m=0,1,2,\ldots,
\end{equation}
and $J_{m+\frac12}(\lambda)$ admits the following explicit representation:
\begin{multline}
\label{half order bessel}
J_{m+\frac12}(\lambda)=\frac1{\sqrt{2\pi}}\sum_{n=1}^{m+1}\beta_n^m\left[\frac{e^{i\lambda}}{i^{n-m}\lambda^{n-\frac12}}+(-1)^{n+m}\frac{e^{-i\lambda}}{i^{n-m}\lambda^{n-\frac12}}\right],\\
\lambda\in\mathbb C\smallsetminus\{0\},\quad m=0,1,2,\ldots,
\end{multline}
where the coefficients $\beta_n^m$ are defined by equations \eqref{betas}.
\end{lemma}
\textbf{Proof} Recall that the finite Fourier transforms of the Legendre polynomials $\hat P_m(\lambda)$ can be expressed in the form
\begin{equation}
\label{5c}
\hat P_m(\lambda)=\frac1{i^m}\sqrt{\frac{2\pi}\lambda}J_{m+\frac12}(\lambda).
\end{equation}
This equation, together with the results of Theorem \ref{Legendre}, immediately gives equation \eqref{half order bessel}.$\hfill\square$

\section{The modified Helmholtz equation in the interior of a square}

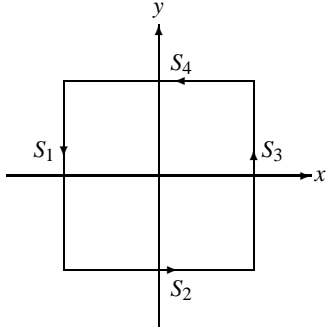
\begin{figure}
\begin{center}
\setlength{\unitlength}{0.5cm}
\begin{picture}(10,10)(-5,-5)
\put(-4,0){\line(1,0){8}}
\put(0,-4){\line(0,1){8}}
\put(-2.5,2.5){\line(1,0){5}}
\put(-2.5,-2.5){\line(1,0){5}}
\put(2.5,-2.5){\line(0,1){5}}
\put(-2.5,-2.5){\line(0,1){5}}
\put(-3.3,0.5){\small$S_1$}
\put(2.7,0.5){\small$S_3$}
\put(0.3,-3.2){\small$S_2$}
\put(0.3,2.8){\small$S_4$}
\put(-2.5,0.6){\vector(0,-1){0.1}}
\put(2.5,0.6){\vector(0,1){0.1}}
\put(0.4,-2.5){\vector(1,0){0.1}}
\put(0.5,2.5){\vector(-1,0){0.1}}
\put(0,3.5){\vector(0,1){0.5}}
\put(3.5,0){\vector(1,0){0.5}}
\put(4.1,-0.15){\small$x$}
\put(-0.15,4.3){\small$y$}
\end{picture}
\caption{The square with sides of length 2.\label{1}}
\end{center}
\end{figure}

Let $u(x,y)$ satisfy the modified Helmholtz equation
\begin{equation}
\label{mH}
\frac{\partial^2u}{\partial x^2}+\frac{\partial^2u}{\partial y^2}-4u=0,
\end{equation}
in the interior of a convex polygon. Associated with each side $\{S_j\}_1^n$ of this $n$-gon, define the following functions
\begin{multline}
\hat u_j(\lambda)=\int_{S_j}\!e^{-i\left(\lambda z-\frac{\bar z}\lambda\right)}\!\left\{\!\left[-u_y+\!\left(\lambda+\frac1\lambda\right)u\right]\!dx+\left[u_x+\left(i\lambda+\frac1{i\lambda}\right)u\right]\!dy\right\}\!,\\
\quad j=1,2,\ldots,n,\quad\lambda\in\mathbb C\smallsetminus\{0\},
\end{multline}
and
\begin{multline}
\tilde u_j(\lambda)=\int_{S_j}\!e^{i\left(\lambda\bar z-\frac z\lambda\right)}\!\left\{\!\left[-u_y+\!\left(\lambda+\frac1\lambda\right)u\right]\!dx+\left[u_x-\left(i\lambda+\frac1{i\lambda}\right)u\right]\!dy\right\}\!,\\
\quad j=1,2,\ldots,n,\quad\lambda\in\mathbb C\smallsetminus\{0\}.
\end{multline}
Then, the functions $\{\hat u_j(\lambda)\}_1^n$ and $\{\tilde u_j(\lambda)\}_1^n$ satisfy the global relations \cite{SSF}
\begin{equation}
\label{GRs}
\sum_{j=1}^n\hat u_j(\lambda)=0,\quad\sum_{j=1}^n\tilde u_j(\lambda)=0,\quad\lambda\in\mathbb C\smallsetminus\{0\}.
\end{equation}
For brevity of presentation, we consider the simplest possible polygon, namely the square with corners at (Figure \ref{1})
\[
(-1,1),\quad(-1,1),\quad(1,-1),\quad(1,1).
\]
For the sides $S_1$, $S_2$, $S_3$, $S_4$, we have respectively
\begin{equation}
z=-1+iy,\quad z=x-i,\quad z=1+iy,\quad z=x+i.
\end{equation}
Hence, taking into account the orientations of the sides, we find the following expressions:
\begin{subequations}
\begin{alignat}1
&\hat u_1(\lambda)=e^{\left(i\lambda+\frac1{i\lambda}\right)}\int_{+1}^{-1}e^{\left(\lambda+\frac1\lambda\right)y}\left[u_x^{(1)}+\left(i\lambda+\frac1{i\lambda}\right)u^{(1)}\right]dy,\\
&\hat u_2(\lambda)=e^{-\left(\lambda+\frac1\lambda\right)}\int_{-1}^{+1}e^{\left(-i\lambda-\frac1{i\lambda}\right)x}\left[-u_y^{(2)}+\left(\lambda+\frac1\lambda\right)u^{(2)}\right]dx,\\
&\hat u_3(\lambda)=e^{-\left(i\lambda+\frac1{i\lambda}\right)}\int_{-1}^{+1}e^{\left(\lambda+\frac1\lambda\right)y}\left[u_x^{(3)}+\left(i\lambda+\frac1{i\lambda}\right)u^{(3)}\right]dy,\\
&\hat u_4(\lambda)=e^{\left(\lambda+\frac1\lambda\right)}\int_{+1}^{-1}e^{\left(-i\lambda-\frac1{i\lambda}\right)x}\left[-u_y^{(4)}+\left(\lambda+\frac1\lambda\right)u^{(4)}\right]dx.
\end{alignat}
\end{subequations}
Let $\hat D_j$ and $\hat N_j$ denote the parts of $\hat u_j$ corresponding to the Dirichlet and Neumann boundary values. Then,
\begin{subequations}
\begin{alignat}1
&\hat u_1(\lambda)=-e^{\left(i\lambda+\frac1{i\lambda}\right)}\hat N_1(\lambda)-\left(i\lambda+\frac1{i\lambda}\right)e^{\left(i\lambda+\frac1{i\lambda}\right)}\hat D_1(\lambda),\\
&\hat u_2(\lambda)=-e^{-\left(\lambda+\frac1\lambda\right)}\hat N_2(-i\lambda)+\left(\lambda+\frac1\lambda\right)e^{-\left(\lambda+\frac1\lambda\right)}\hat D_2(-i\lambda),\\
&\hat u_3(\lambda)=e^{-\left(i\lambda+\frac1{i\lambda}\right)}\hat N_3(\lambda)+\left(i\lambda+\frac1{i\lambda}\right)e^{-\left(i\lambda+\frac1{i\lambda}\right)}\hat D_3(\lambda),\\
&\hat u_4(\lambda)=e^{\left(\lambda+\frac1\lambda\right)}\hat N_4(-i\lambda)-\left(\lambda+\frac1\lambda\right)e^{\left(\lambda+\frac1\lambda\right)}\hat D_4(-i\lambda).
\end{alignat}
\end{subequations}
For simplicity, we consider the symmetric Dirichlet boundary value problem:
\begin{subequations}
\label{bvp}
\begin{alignat}4
&u^{(1)}&&\!=u(-1,y)&&\!=\cosh(1)\cosh(\sqrt{3}y)\!+\!\cosh(\sqrt{3})\cosh(y),\, &&-1\!<\!y\!<\!1;\\
&u^{(3)}&&\!=u(1,y) &&\!=\cosh(1)\cosh(\sqrt{3}y)\!+\!\cosh(\sqrt{3})\cosh(y),\, &&-1\!<\!y\!<\!1;\\
&u^{(2)}&&\!=u(x,-1)&&\!=\cosh(1)\cosh(\sqrt{3}x)\!+\!\cosh(\sqrt{3})\cosh(x),\, &&-1\!<\!x\!<\!1;\\
&u^{(4)}&&\!=u(x,1) &&\!=\cosh(1)\cosh(\sqrt{3}x)\!+\!\cosh(\sqrt{3})\cosh(x),\, &&-1\!<\!x\!<\!1.
\end{alignat}
\end{subequations}
Then,
\begin{equation}
u(x,y)=u(-x,y),\quad u(x,y)=u(x,-y),\quad u(x,y)=u(y,x).
\end{equation}
Thus,
\begin{equation}
u_x^{(3)}=-u_x^{(1)},\quad u_y^{(4)}=-u_y^{(2)},\quad \left.u_y^{(2)}=u_x^{(1)}\right|_{(x,y)\leftrightarrow(y,x)}.
\end{equation}
Hence, the first of the global relations \eqref{GRs} becomes
\begin{subequations}
\label{GR}
\begin{multline}
\label{GR 1}
\cos\left(\!\lambda-\frac1\lambda\!\right)\hat N_1(\lambda)+\cos\left(\!i\lambda-\frac1{i\lambda}\!\right)\hat N_1(-i\lambda)\\
\quad= \left(\!\lambda-\frac1\lambda\!\right)\sin\left(\!\lambda-\frac1\lambda\!\right)\hat D_1(\lambda)
+\left(\!i\lambda-\frac1{i\lambda}\!\right)\sin\left(\!i\lambda-\frac1{i\lambda}\!\right)\hat D_1(-i\lambda),\\
\lambda\in\mathbb{C}\smallsetminus\{0\}.
\end{multline}
\newpage
\noindent The simplest way to obtain the second of the global relations \eqref{GRs} is to take the Schwartz conjugate of equation \eqref{GR 1} (i.e. to take the complex conjugate of \eqref{GR 1} and then to replace $\bar\lambda$ with $\lambda$). This yields the equation
\begin{multline}
\label{GR 2}
\cos\left(\!\lambda-\frac1\lambda\!\right)\hat N_1(\lambda)+\cos\left(\!i\lambda-\frac1{i\lambda}\!\right)\hat N_1(i\lambda)\\
\quad= \left(\!\lambda-\frac1\lambda\!\right)\sin\left(\!\lambda-\frac1\lambda\!\right)\hat D_1(\lambda)
+\left(\!i\lambda-\frac1{i\lambda}\!\right)\sin\left(\!i\lambda-\frac1{i\lambda}\!\right)\hat D_1(i\lambda),\\
\lambda\in\mathbb{C}\smallsetminus\{0\}.
\end{multline}
\end{subequations}

Using $N$ basis functions to approximate $u_x^{(1)}$, equations \eqref{GR} yield 2 equations for $N$ unknowns.

Regarding the numerical solution of the global relations \eqref{GR}, Fourier basis functions \cite{FFSS}, \cite{FFX}--\cite{SSF}, as well as Chebyshev and Legendre polynomials \cite{D}, \cite{FF}, \cite{SFFS}, have been used to approximate $u_x^{(1)}$. In most of the earlier papers the collocation points $\lambda\in\mathbb{C}\smallsetminus\{0\}$ were chosen to lie on the rays in the complex $\lambda$-plane which are parallel to the edges of the polygon and its reflection in the imaginary axis. Recently B Fornberg and collaborators introduced the use of the so-called Halton nodes, \cite{D}, \cite {FF}.

It appears that the most efficient numerical method involves (a) approximating the data $u_x^{(1)}$ in terms of Legendre polynomials (following Fornberg) and (b) using the collocation points employed in our earlier work \cite{SSF} where ideas of A G Sifalakis and collaborators \cite{SFFS} for the Laplace equation were extended to the modified Helmholtz equation. Furthermore, in order to ensure that the collocation matrix remains well conditioned as the number $N$ of basis functions increases, it is important following Fornberg to (i) divide each row, as well as each column, of the collocation matrix by its $l^1$-norm \cite{D} and (ii) to ``over-determine'' the linear system by choosing the number of collocation points to be about the same as the number of unknowns.

Numerical experiments suggest that Legendre polynomials yield spectral accuracy rather than the algebraic accuracy found with a Fourier basis. Furthermore, as a result of choosing the collocation points to be on the above rays, the semi-block circulant structure of the collocation matrix for regular polygons, demonstrated for the Laplace equation in \cite{SSP}, is preserved in modified Helmholtz equation as well.

\begin{figure}
\[
\begin{array}{ccc}
\!\!\!\!\!\!\!\!\!\!\!\!\epsfig{figure=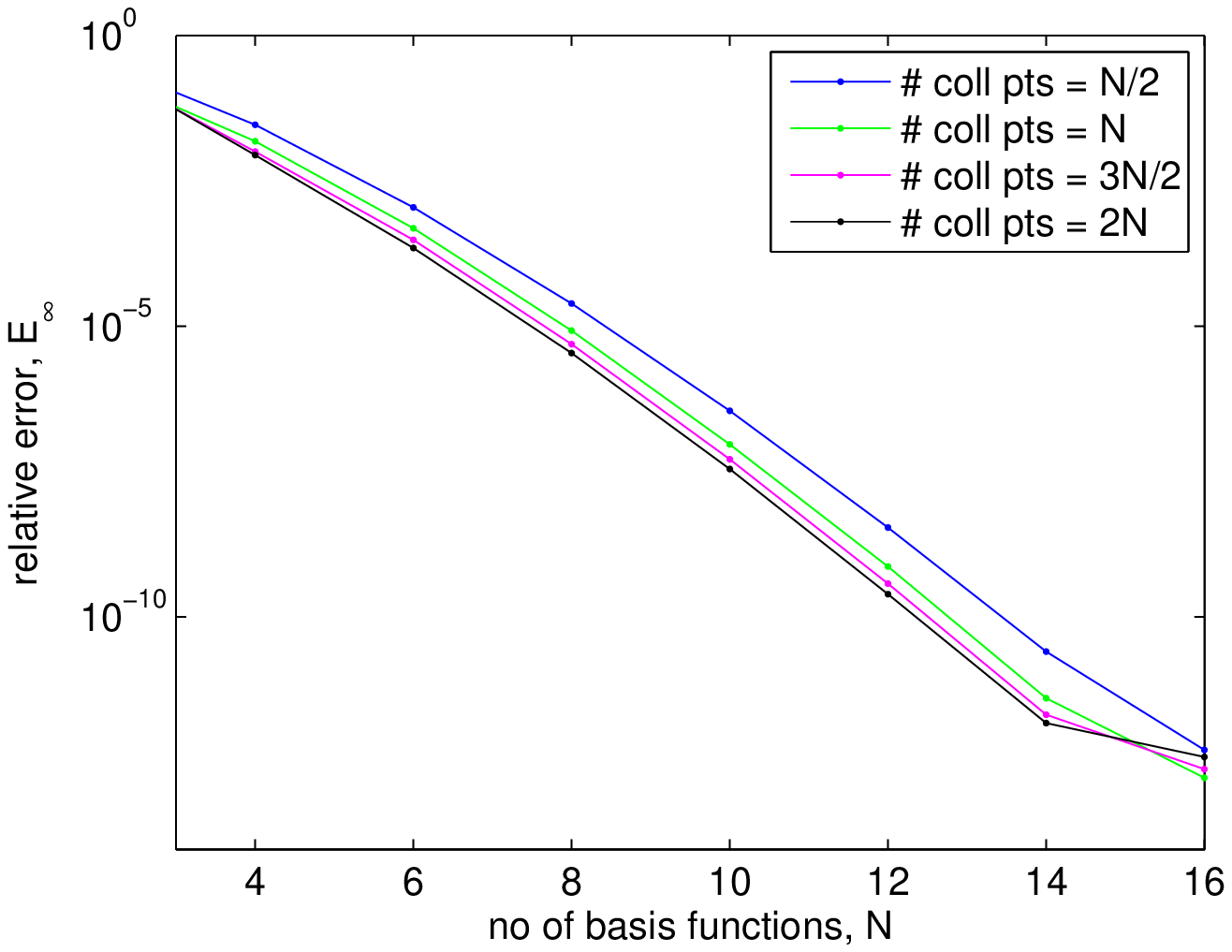,width=5.8cm}&\!\!\!\!\!\!\!\!\!\!\!\!\!\!\!\!\!\!\!\!\!\!\!\!&\epsfig{figure=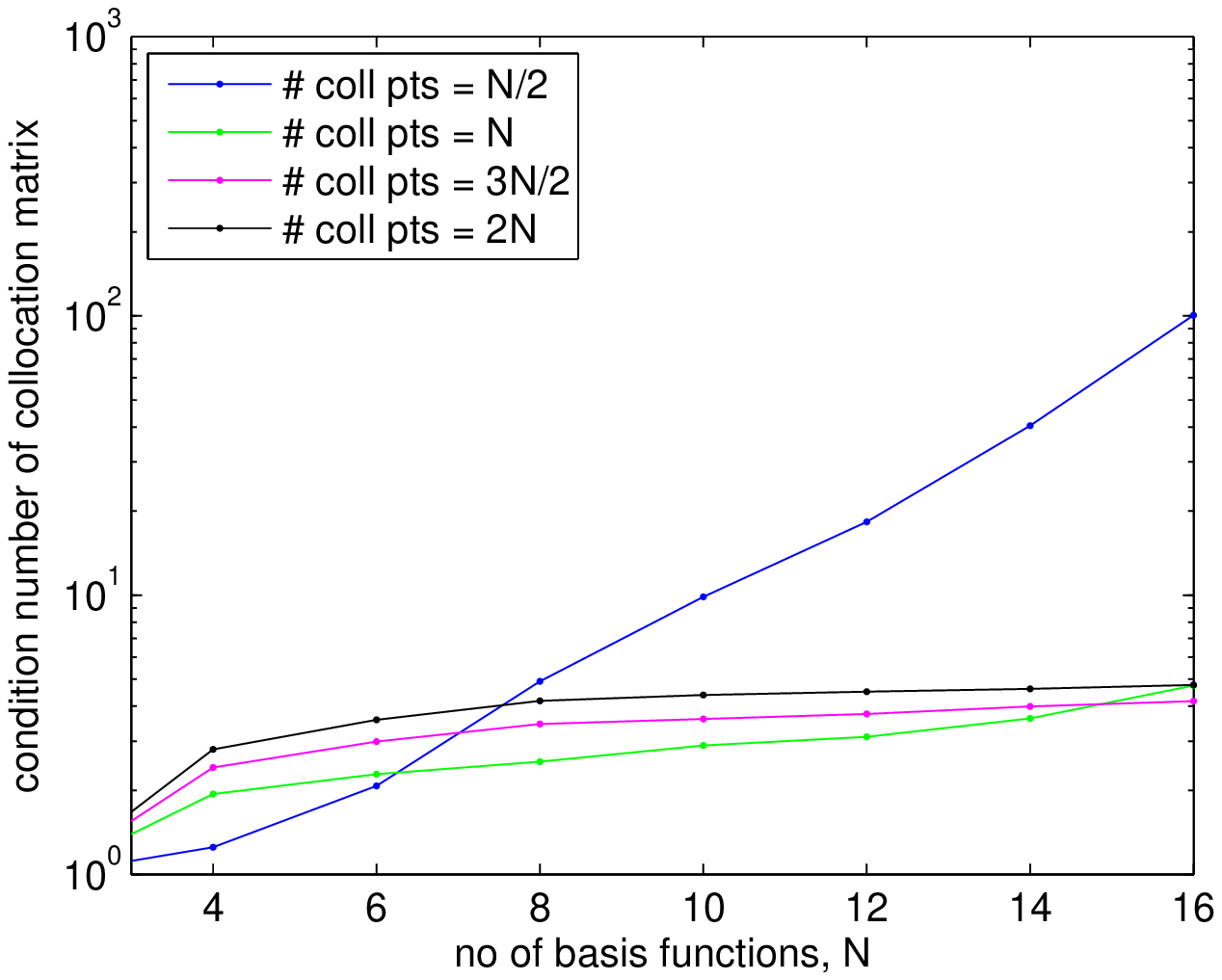,width=5.8cm}\!\!\!\!\!\!\!\!\!\!\!\!\!\!\!\!\!\!\!\!\!\!\!\!
\end{array}
\]
\caption{Numerical solution of the global relations \eqref{GR} for the symmetric Dirichlet boundary value problem \eqref{bvp}.\label{2}}
\end{figure}

Plots of the relative error $E_\infty$ (defined in \cite{SSF}), as well as of the matrix condition number as a function of $N$, for $N/2$, $N$, $3N/2$ and $2N$ collocation points, are presented in Figure \ref{2}. The rectangular collocation matrix was inverted by using the ``backslash'' command in Matlab. It is clear that over-determining the linear system by a factor of 2 is sufficient to achieve very good matrix conditioning.

\section{Acknowledgements}
ASF is grateful to EPSRC, UK and to Onassis foundation, USA for partial support. SAS wishes to thank \mbox{Prof B Fornberg} for many useful discussions.

\end{document}